\documentclass[reqno]{amsart}
\usepackage{microtype}

\usepackage[shortlabels]{enumitem}
\setlist[enumerate]{label={\arabic*.}}

\usepackage{amssymb}
\usepackage{bm}  
\usepackage[dvipsnames]{xcolor}
\usepackage
[colorlinks=true,linkcolor=Maroon,citecolor=OliveGreen]
{hyperref}
\usepackage{bookmark}
\usepackage[abbrev,shortalphabetic]{amsrefs}  
\usepackage{cleveref} 
\usepackage{tikz}

\usepackage{xcolor}

\usepackage[norefs, nocites]{refcheck}

\makeatletter
\newcommand{\refcheckize}[1]{%
  \expandafter\let\csname @@\string#1\endcsname#1%
  \expandafter\DeclareRobustCommand\csname relax\string#1\endcsname[1]{%
    \csname @@\string#1\endcsname{##1}\wrtusdrf{##1}}%
  \expandafter\let\expandafter#1\csname relax\string#1\endcsname
}
\makeatother

\refcheckize{\cref}
\refcheckize{\Cref}

\numberwithin{equation}{section}

\newtheorem{theorem}{Theorem}
\newtheorem*{theorem*}{Theorem}
\newtheorem{lemma}[theorem]{Lemma}

\theoremstyle{definition}

\renewcommand{\phi}{\varphi}

\renewcommand{\subset}{\subseteq}

\newcommand\opr[1]{\operatorname{#1}}

\def\PGL{\opr{PGL}}
\def\PSL{\opr{PSL}}

\def\P{\mathbf{P}}

\newcommand\mat[4]{\begin{pmatrix} #1 & #2 \\ #3 & #4 \end{pmatrix}}

\usepackage{todonotes}

\begin{document}
    
    \title{Sharply transitive sets in $\PGL_2(K)$}
    
    \author{Sean Eberhard}
    \address{Sean Eberhard, Centre for Mathematical Sciences, Wilberforce Road, Cambridge CB3~0WB, UK}
    \email{eberhard@maths.cam.ac.uk}
    
    \thanks{SE has received funding from the European Research Council (ERC) under the European Union’s Horizon 2020 research and innovation programme (grant agreement No. 803711).
    }
    
    \begin{abstract}
        Here is a simplified proof that every sharply transitive subset of $\PGL_2(K)$ is a coset of a subgroup.
    \end{abstract}
    
    \maketitle
    
    Let $G$ be a group acting on a set $\Omega$ on the left.
    A subset $S \subset G$ is called \emph{sharply transitive}, or \emph{regular},
    if for every $x, y \in \Omega$ there is a unique $g \in G$ such that $gx = y$.
    For example, sharply transitive subsets of $S_n$ can be identified with $n \times n$ Latin squares,
    while sharply 2-transitive subsets, i.e., sets which are sharply transitive for the action on ordered pairs,
    can be identified with affine planes of order $n$.
    
    In this note we consider $G = \PGL_2(K)$ and its action on the projective line $\P^1 = \P^1(K)$,
    where $K$ is a finite field of order $q$.
    It follows from Dickson's classification of the subgroups of $\PSL_2(K)$ that the only regular subgroups of $\PGL_2(K)$ are
    \begin{enumerate}
        \item cyclic groups $C_{q+1}$,
        \item dihedral groups $D_{(q+1)/2}$ ($q$ odd),
        \item $A_4$ ($q = 11$),
        \item $S_4$ ($q = 23$),
        \item $A_5$ ($q = 59$)
    \end{enumerate}
    (see \cite{bonisoli} or \cite{VM}),
    and there is a single conjugacy class of subgroups in each case.
    Remarkably, other than these subgroups and their cosets, there are no further regular subsets of $G$.
    
    \begin{theorem}
        [\cites{bader-lunardon, thas, bonisoli-korchmaros}]
        \label{thm:main}
        If $S \subset \PGL_2(K)$ is sharply transitive on $\P^1$ and $1 \in S$ then $S$ is a subgroup.
    \end{theorem}
    
    This result was originally conjectured by Bonisoli in \cite{bonisoli}.
    The complete proof is spread across several papers:
    in \cites{bader-lunardon, thas}, Bader, Lundardon, and Thas classified flocks of the hyperbolic quadric in $\P^3$,
    and the equivalence with regular sets was noted in \cite{bonisoli-korchmaros}.
    The original proof is somewhat involved.
    See \cite{thas2} for a more recent summary.
    A partly simplified proof is given by Durante and Siciliano~\cite{durante-siciliano}, but for the main technical step the reader is referred to \cite{thas}.
    Here we give a short, self-contained, direct proof avoiding the detour through flocks.
    Still, the interested reader will find close analogies with \cites{durante-siciliano, thas}.
    
    The key step is the following lemma, which can be seen as a version of Segre's ``lemma of the tangents'' (see \cite{segre}*{(2)}).
    
    \begin{lemma}
    Suppose $S$ is a regular subset of $\PGL_2(K)$ containing the elements
    \begin{align*}
        &g_1 = \mat{0}{b_1}{c_1}{d_1},
        &&g_2 = \mat{a_2}{0}{c_2}{d_2},\\
        &g_3 = \mat{a_3}{b_3}{0}{d_3},
        &&g_4 = \mat{a_4}{b_4}{c_4}{0}.
    \end{align*}
    Then
    \[
        b_1 d_2 a_3 c_4 = c_1 a_2 d_3 b_4.
    \]
    \end{lemma}
    
    Here and below elements of $\PGL_2(K)$ are written as matrices,
    understanding that the expression is determined only up to a scalar.
    Elements of $\P^1$ will be written as $(x : y)$.
    
    \begin{proof}
    Assume first that $g_1, g_2, g_3, g_4$ are distinct. Label the elements of $S$ as
    \[
        \mat{a_i}{b_i}{c_i}{d_i} \qquad (i = 1, \dots, q+1).
    \]
    By regularity, the columns $(a_i : c_i)$ and $(b_i : d_i)$ both trace out $\P^1$.
    Since $S^{-1}$ is also regular, the same is true of the rows $(a_i : b_i)$ and $(c_i : d_i)$.
    By excluding $0$ and $\infty$ in each case, it follows that the products
    \[
        \prod_{i \neq 1, 2} a_i / b_i,
        \prod_{i \neq 2, 4} b_i / d_i,
        \prod_{i \neq 3, 4} d_i / c_i,
        \prod_{i \neq 1, 3} c_i / a_i
    \]
    are each equal to the product of all the nonzero elements of $K$, which is $-1$.
    By multiplying these together we get
    \[
        \frac{b_1 d_1}{d_1 c_1}
        \frac{d_2 c_2}{c_2 a_2}
        \frac{a_3 b_3}{b_3 d_3}
        \frac{c_4 a_4}{a_4 b_4}
        \prod_{i > 4} \frac{a_i}{b_i} \frac{b_i}{d_i} \frac{d_i}{c_i} \frac{c_i}{a_i} = 1.
    \]
    Simplifying gives the claimed relation.
    
    The only possible coincidences between $g_1, g_2, g_3, g_4$ are $g_1 = g_4$ and $g_2 = g_3$.
    In these cases the proof is similar.
    \end{proof}
    
    We can now prove \Cref{thm:main}.
    Let $S \subset \PGL_2(K)$ be a regular set containing $1$.
    Let $g, h \in S \setminus \{1\}$.
    It suffices to prove $gh \in S$.
    By regularity, for each $x \in \P^1$ there is a unique $k \in S$ such that
    \[
        kh^{-1}x = gx.
    \]
    Moreover $k \neq h$, because $g$ has no fixed points (likewise $k \neq g$).
    Since $\P^1$ has size $q+1$ and there are at most $q$ possibilities for $k$,
    there is some $k \in S$ such that $g^{-1}kh^{-1}$ has at least two fixed points, say $x$ and $y$.
    Since $x$ and $gx$ are distinct, we may coordinatize $\P^1$ in such a way that
    \[
        x = (1:0),\quad
        gx = (0:1).
    \]
    Since $h$ has no fixed points, $h (1:0) \neq (1:0)$, so there is some $u \in \PGL_2(K)$ of the form
    \[
        u = \mat1*01
    \]
    such that $hu(0:1) = (1:0)$.
    For this $u$ we have
    \[
        ku(0:1) = kh^{-1}(1:0) = g(1:0) = (0:1).
    \]
    Thus $Su$ is a regular set containing the elements
    \begin{align*}
        &gu = \mat0b1* ,
        &ku = \mat{a}0*1,\\
        &u = \mat1*01,
        &hu = \mat*1c0,
    \end{align*}
    for unknown nonzero values $a, b, c$ (and further unknown values hidden by $*$).
    By the lemma, we must have $a = bc$.
    Hence
    \begin{align*}
        hk^{-1}g
        &= (hu) (ku)^{-1} (gu) u^{-1} \\
        &= \mat*1c0 \mat10*a \mat0b1* \mat1*01 \\
        &= \mat{a}*0{bc} \\
        &= \mat1*01.
    \end{align*}
    Since $hk^{-1}g$ has at least two fixed points ($x$ and $y$), we must have $hk^{-1}g=1$.
    Hence $gh = k \in S$, as required.
    
    \bibliography{refs}
\end{document}